%========================================================================
% A TeX file for the paper "The Locally Fine Coreflection and Normal
%             Covers in the Products of Partition-complete Spaces"
%
% You must have a text file named "part.ref" in the same directory
% with this file. (Simply create an empty file.) Also, you need to 
% compile this file twice in order to get the reference numbers defined. 
\newif\ifmakeref 
\makereffalse    % set this to true for updated references
%\makereftrue
%
% Aarno Hohti & Miroslav Husek & Jan Pelant started  Dec 1997
% Last modification A.Hohti Jan 2002
% 
% The modifications to comply with the referee for Topology & Applications
% are marked by "%%!!"
%=======================================================================
%
%\magnification = \magstep1     % for weak eyes!
%
\baselineskip = 15.5pt plus .5pt minus .5pt   % for improved readability!
%
%========================================================================
\font\bigbf=cmbx10   scaled \magstep1     % for titles
                       % for file contents inside text
                         % for group names
\font\smallcaps = cmcsc10                 % for headlines
\font\bs=cmssbx10                         % bold sans serif
%===============================================================
% This is a TeX macro for switching to nine point fonts. Taken
% from  the TeXbook, pages 413 - 415. (For special cases only)
%===============================================================
%
\ifx\fonts\cmfonts
\font\ninerm=cmr9
\font\ninei=cmmi9
\font\ninesy=cmsy9
\font\ninebf=cmbx9
\font\ninett=cmtt9
\font\nineit=cmti9
\font\ninesl=cmsl9
\else
\font\ninerm=amr9
\font\ninei=ammi9
\font\ninesy=amsy9
\font\ninebf=ambx9
\font\ninett=amtt9
\font\nineit=amti9
\font\ninesl=amsl9
\fi
\skewchar\ninei='177
\skewchar\ninesy='60
\skewchar\ninett=-1
\newskip\tglue
\def\ninepoint{\def\rm{\fam0\ninerm}
       \textfont0=\ninerm
       \textfont1=\ninei
       \textfont2=\ninesy
       \textfont\itfam=\nineit \def\it{\fam\itfam\nineit}%
       \textfont\slfam=\ninesl \def\sl{\fam\slfam\ninesl}%
       \textfont\ttfam=\ninett \def\tt{\fam\ttfam\ninett}%
       \textfont\bffam=\ninebf \def\bf{\fam\bffam\ninebf}%
       \tt\tglue=.5em plus.25em minus .15em
       \normalbaselineskip=11pt
       \setbox\strutbox=\hbox{\vrule height8pt depth3pt width0pt}%
       \let\sc=\sevenrm \let\big=\ninebig \normalbaselines\rm}%
                           % abstract, references
\input amssym
\input xypic
\input part.ref
\def\fonts{cmfonts}
\advance\tolerance by 100

%\AdobeTimesRoman 
%========================================================================
%
%========================================================================
\def\runningtitlestring{Locally Fine
                    Coreflection in Product Spaces}
% stylistic definitions
\newcount\sectioncount
\sectioncount = 0
\def\ThisSection{\number\sectioncount}
\def\section#1{\vskip0pt plus .1\vsize
    \penalty-250\vskip0pt plus-.1\vsize\bigskip
    \global\advance\sectioncount by 1
    \noindent{\bf \number\sectioncount. #1}\nobreak\message{#1}}
\def\SectionBreak{\vskip0pt plus .1\vsize
    \penalty-250\vskip0pt plus-.1\vsize\null\bigskip}

\def\abstract#1{{\ninepoint\bigskip\centerline{\hbox{
       \vbox{\hsize=4.85truein{\noindent ABSTRACT.\enspace{#1}}
            }}}}}
\def\keywords#1{{\ninepoint\bigskip\centerline{\hbox{
       \vbox{\hsize=5.75truein{\noindent KEYWORDS.\enspace#1}
            }}}}}
\def\AmsClass#1#2{{\ninepoint\bigskip\centerline{\hbox{
       \vbox{\hsize=5.75truein{\noindent
       AMS #1 Subject Classification:\enspace#2}
            }}}}}
\def\th #1 #2: #3\par{\medbreak{\bf#1 #2:
\enspace}{\sl#3\par}\par\medbreak}
\def\co #1 #2: #3\par{\medbreak{\bf#1 #2:
\enspace}{\sl#3\par}\par\medbreak}
\def\le #1 #2: #3\par{\medbreak{\bf #1 #2:
\enspace}{\sl #3\par}\par\medbreak}
\def\rem #1 #2. #3\par{\medbreak{\bf #1 #2.
\enspace}{#3}\par\medbreak}
\def\proof{{\bf Proof}.\enspace}

\def\sqr#1#2{{\vcenter{\hrule height.#2pt
      \hbox{\vrule width.#2pt height#1pt \kern#1pt
         \vrule width.#2pt}
       \hrule height.#2pt}}}
\def\eop{\mathchoice\sqr34\sqr34\sqr{2.1}3\sqr{1.5}3}
\overfullrule=0pt
% this definition is from the TeXbook, Ex. 21.3 :
\def\boxit#1{\vbox{\hrule \hbox{\vrule \kern2pt
                 \vbox{\kern2pt#1\kern2pt}\kern2pt\vrule}\hrule}}
% A macro for making references and blocks.
%
\newdimen\refindent\newdimen\plusindent
\newdimen\refskip\newdimen\tempindent
\newdimen\extraindent
\newcount\refcount
\newwrite\reffile
\def\beginref{\ifmakeref\immediate\openout\reffile = part.ref\else\fi}
\def\endref{\ifmakeref\immediate\closeout\reffile\else\fi}
\def\ifundefined#1{\expandafter\ifx\csname#1\endcsname\relax}
\def\referto[#1]{\ifundefined{#1}[?]\else[\csname#1\endcsname]\fi}
\def\refered[#1]{\ifundefined{#1}[?]\else\csname#1\endcsname\fi}
%
% \refskip has to be set by the user! Otherwise \parindent is
% used, in accordance with \item.
%
\refcount=0
\def\ref#1:#2.-#3[#4]{\ninepoint % switch to nine point fonts
\advance\refcount by 1
\setbox0=\hbox{[\number\refcount]}\refindent=\wd0
\plusindent=\refskip\extraindent=\refskip
\advance\plusindent by -\refindent\tempindent=\parindent %
\parindent=0pt\par\hangindent\extraindent %
 [\number\refcount]\hskip\plusindent #1:{\sl#2},#3
\parindent=\tempindent
\ifmakeref
\immediate\write\reffile{\string\expandafter\def\noexpand\csname#4\endcsname%
{\number\refcount}}\else\fi}
\refskip=\parindent
\def\markplace#1#2{\immediate\write\reffile{\string\expandafter\def\noexpand%
\csname#1\endcsname{#2}}}

\def\makeknown#1#2{\expandafter\gdef
            \csname#2\endcsname{\hbox{\csname#1\endcsname #2}}}
%%%%%%%%%%%%%%%%%%%%%%%%%%%%%%%%%%%%
% globally known `text constants'
\makeknown{rm}{Cov}
\makeknown{rm}{int}
\makeknown{rm}{Ord}
\makeknown{rm}{Rank}
\makeknown{rm}{Root}
\makeknown{rm}{End}
\makeknown{rm}{Level}
\makeknown{rm}{Seg}
\makeknown{rm}{pr}
\makeknown{rm}{top}
\def\c#1{{\cal #1}}
\def\st{\hbox{\rm St}}

\def\KK{\cal K}

%

%--------------------------------------------------------
\def\specialheadlines{
        \headline={\vbox{\line{
                     \ifnum\pageno<2 \ffolio
                     \else\rightheadline
                     \fi}\bigskip
                     }}}
        \def\ffolio{\hfil}
        \def\rightheadline{\hfil
                    {\smallcaps\runningtitlestring}\hfil
                    \hbox{\rm\folio}}
\nopagenumbers
\specialheadlines
\beginref
%--------------------------------------------------------
%
%
%========================================================================
                              % HEADER
\null
\vskip1.5truecm
  \centerline{\bigbf The Locally Fine Coreflection and Normal Covers}
  \medskip
  \centerline{\bigbf in the Products of Partition-complete Spaces}
\vskip.75truecm
 \centerline{\sl A.\ Hohti, University of Helsinki}
 \centerline{\sl M.\ Hu\v sek, Charles University, Prague}
 \centerline{\sl J.\ Pelant, Czech Academy of Sciences\footnote{*}
   {\ninepoint The first author is grateful to the {\it Finnish Academy of Sciences}
    for a 1997 grant which partially supported this research, and 
    the second and third authors are grateful to the {\it Grant
    Agency of the Czech Republic} (GA \v CR 201/00/1466).}
}
\vskip.5truecm
%            --------------------------------------------------------
 %           --------------------------------------------------------
\vskip.5truecm
 %        --------------------------------------------------------------
\abstract{We prove that the countable product of supercomplete spaces having
          a countable closed cover consisting of partition-complete
          subspaces is supercomplete with respect to its metric-fine
          coreflection. Thus, countable products of
          $\sigma$-partition-complete paracompact spaces are again paracompact. On the
          other hand, we show (\refered[ThNormalCovers]) that in arbitrary products of
          partition-complete paracompact spaces, all
          %%!!$\sigma$-partition-complete paracompact spaces, all
          normal covers can be obtained via the locally fine
          coreflection of the product of fine uniformities. 
          %% THE FOLLOWING NOT YET IN THIS VERSION!
          %For 1st-countable factors one obtains new extension theorems
          %of continuous functions defined on the subspaces of the
          %product. 
          These results extend those given in
          \referto[Alstera],
          \referto[Friedlera],
          \referto[Frolika],
          \referto[HohtiPelantb],
          \referto[Pelanta],
          %\referto[HusekPelant],
          \referto[HohtiYun],
          \referto[Pleweb].
          }

%
%\smallskip
\AmsClass{2000}{Primary 54B10, 54E15, Secondary 54D20.}
\vskip.35truecm
%========================================================================
%
\section{Introduction.} It was shown by the third author in
\referto[Pelanta] that {\it locally fine spaces} are subspaces of fine
spaces. This result was obtained by embedding an arbitrary
uniform space into the product of complete metric spaces and
by showing that the {\it locally fine coreflection}
\referto[Ginsburga] of such a product
uniformity contains all normal covers of the product. The important
insight that emerged from \referto[Pelanta] was the connection of
{\it Noetherian trees} with the covers of locally fine spaces. The
locally fine coreflection $\lambda\mu$ of a given {\it pre-uniformity}
$\mu$ (i.\ e., a filter of coverings) may be described as an attempt
to obtain the topology of a given space by means of
iterated combinatorial refinements (generalized subdivisions)
of coverings, starting from $\mu$. %%!! as the given base. 
The covers obtained in this way are recursively related to $\mu$ %%!!the base
in the sense that for each $\c{U}\in \lambda\mu$ there is a
Noetherian tree $T$ of subsets of $X$, obtained through
iterated applications of uniform refinements to the elements
of the tree, such that the {\it maximal elements} of $T$ form a
refinement of $\c{U}$. For {\it supercomplete spaces} \referto[Isbella],
and in particular for complete metric spaces, this procedure
reaches {\it all open covers} of the space. The complexity of
such refinements was studied in \referto[Hohtic] and
\referto[HohtiPelanta]. Related methods were used by the
second and third authors to obtain extension theorems for
continuous functions defined on subsets of products of metrizable
spaces \referto[HusekPelant].

On the
other hand, the condition that the locally fine coreflection
of a product contain all open covers was studied by the
first author in a series of papers on supercompleteness
(see, e.\ g., \referto[Hohtib], \referto[Hohtid],
\referto[Hohtif]). For spaces satisfying a structural  %%!! 
recursive condition (C-scattered, \v Cech-scattered)
Noetherian trees were applied to show that paracompactness
(among others) is preserved in countable products of such
spaces. These results extended the initial theorems of
Frol\'\i k \referto[Frolika] and Arhangel'skii
\referto[Arhangelskib] that paracompactness is
preserved in countable products of locally compact spaces,
and similar results from \referto[Rudina] (scattered spaces)  %%!! by --> from
and \referto[Friedlera] (C-scattered spaces), including the
result of \referto[Alstera] on C-scattered Lindel\"of spaces.
However, and more importantly, the new results showed that the %%!! they --> new results
open covers of such product could be refined by ones
generated recursively from open covers of the factors.

For paracompact products, the above condition on the locally fine %%!!
coreflection is equivalent to the {\it spatiality}
of the corresponding product of {\it locales}, noted         
by Isbell in \referto[Isbellc]. Spatiality was extensively studied %%!!
by Plewe in \referto[Plewea], who used the productivity
of paracompactness in locales and spatiality to
extend (\referto[Pleweb]) the above topological
results to countable products of $C_\delta$-absolute spaces. %%!!
These spaces are equivalent to those with a {\it complete    %%!!
exhaustive sieve} studied by Michael in \referto[Michaelb]
and to the {\it partition-complete} spaces
of Telg\'arsky and Wicke
in \referto[TelgarskyWicke].
Unfortunately, spatiality does not extend to
spaces which are merely countable unions
(e.\ g., the rationals) of closed partition-complete parts.
As an extreme case, $\Bbb Q\times \Bbb Q$ {\it is not} spatial.

It has been recently shown by the first author (\referto[Hohtig])
that for spaces with pre-uniformities, there is a 
simple relationship between the `covering monoid' of
their localic product and the locally fine coreflection
of the product pre-uniformity. The covers consisting of open rectangles 
of the latter product form a monoid
which naturally maps onto a `generating monoid' of 
the localic product.
It follows that for the products of regular 
spaces with supercomplete   %%!! (super)complete
pre-uniformities, the locally fine coreflection of the
product pre-uniformity is supercomplete provided that the
product is spatial. (We call a pre-uniformity supercomplete %%!! 
if it contains every open cover of the underlying space.)
On the other hand, the locally fine pre-uniformities are
equivalent to the so-called {\it formal spaces}, for which
the essential questions about covers are related to their
effective computability. Thus, we come back to our methods
of recursive constructions by means of Noetherian trees.
As we cannot rely on the spatiality of the products, we will
construct the refinements directly. Moreover, we will handle
{\it arbitrary} (in general non-paracompact) products of
%%!!$\sigma$-partition-complete regular spaces
partition-complete regular spaces
and show that their locally fine covers by regular open
sets belong to (i.\ e., are refined by members of) the
locally fine coreflection of any supercomplete uniformities
of the factors.

Before starting our preliminary section for necessary definitions,
we state our two main results:

\th Theorem A (\refered[ThSigmaPartitionComplete]): 
Let $(X_i: i\in \Bbb N)$
be a countable family of $\sigma$-partition-complete, regular
spaces, and let $\mu_i$ be a supercomplete pre-uniformity
for $X_i$. Then
$\lambda m\Pi\mu_i$ is a supercomplete pre-uniformity for
the topological product $\Pi X_i$.

\th Theorem B (\refered[ThNormalCovers]): Let $(X_i: i\in I)$
be a family of partition-complete, completely regular spaces, and for
each $i\in I$, let $\mu_i$ be a supercomplete uniformity
for $X_i$. Then every normal cover $\c{V}$ of
$X = \Pi_I X_i$ belongs to $\lambda\Pi_I\mu_i$.

\section{Preliminaries.} Here we take up definitions of several
tools needed in the main part of this paper. The locally fine
coreflections of uniformities were introduced in
\referto[Ginsburga], and they were constructed internally
by means of consecutive `derivatives' closing the
given (pre-) uniformity under the combinatorial
condition of local fineness. In the theory of
`formal spaces', this closure is considered
via the condition of {\it composition} (\referto[FourmanGrayson])
or {\it transitivity} (\referto[Sambina]) of the
related covering relation. It may be
regarded as an attempt to obtain the underlying topology
of the space through iterated applications of
uniform subdivisions of covers. We recall the
definition of $\lambda$ from \referto[Hohtib] by means
of slowed-down Ginsburg-Isbell derivatives. Let
$\mu,\nu$ be pre-uniformities on a set $X$. Then
(\referto[Riceb]) $\mu/\nu$ denotes the pre-uniformity
generated by all covers of $X$ of the form
$\{U_i\cap V^i_j\}$, where $\{U_i\}\in \mu$ and
for all $i$, $\{V^i_j\}\in \nu$. Let $\mu^{(0)} = \mu$,
$\mu^{(\alpha+1)} = \mu^{(\alpha)}/\mu$, and
$\mu^{(\beta)} = \bigcup\{\mu^{(\alpha)}: \alpha<\beta\}$
for a limit ordinal $\beta$. The first $\mu^{(\alpha)}$
unchanged under this derivation is denoted by
$\lambda\mu$ and called the locally fine coreflection of
$\mu$. It is the minimal pre-uniformity closed under
the operation $\nu\mapsto \nu/\nu$ finer than $\mu$.
If $\c{V}\in \lambda\mu$, we call $\c{V}$ a
$\lambda$-cover (with respect to $\mu$) or directly
a $\lambda\mu$-uniform cover. For a subset $A\subset X$,
a $\lambda$-neighbourhood (with respect to a given
pre-uniformity $\mu$) is a set $N\supset A$ for which
there is $\c{V}\in \lambda\mu$ such that
$\st(A,\c{V})\subset N$. We note here that if the
$\mu_i$ are pre-uniformities and $\c{V}\in \lambda\Pi\mu_i$,
then $\c{V}$ has a refinement $\c{W}\in \lambda\Pi \mu_i$
consisting of basic rectangles, i.\ e., sets of the
form $B = \cap\{\pi_i^{-1}[B_i]: i\in F\}$, where $F$
is a finite set. This result can be easily proved by
using the inductive definition of the operation $\lambda$.

{\bf Noetherian trees.}
Well-founded trees
in the context of locally fine refinements were first
used in \referto[PelantRice] and in an early version of
\referto[Pelanta]. They were consequently applied in
several papers: \referto[HohtiPelanta],
\referto[HohtiPelantb], \referto[HohtiYun], and in
\referto[Hohtih] in the context of well-founded cubical
triangulations refining open covers of infinite-dimensional cubes.

In this paper,
a tree is a set $T$ with a partial order $\leq$ such that
$T$ has a unique minimal element (root) with respect to
$\leq$ and for each $p\in T$, the set of $\leq$-predecessors
of $p$ is linearly ordered and finite. We call a tree
{\it Noetherian} if each of its linearly ordered parts is finite.
To fix some notation, we denote by $S(p)$ the set of all
immediate successors of $p\in T$.
The symbol $\End(T)$ denotes the set of all maximal
elements of $T$. Thus, if $T$ is Noetherian, then
not only is $\End(T)$ non-empty, but each maximal linear
part (branch) of $T$ meets $\End(T)$ in a unique element.

The Ginsburg-Isbell derivatives are related with Noetherian
trees as follows. Let $\mu$ be a pre-uniformity on a
set $X$. We say that a mapping $\varphi: T\to 2^X$ satisfies
the {\it uniform mapping condition} with respect to
$\mu$ if for each $p\in T$, the family
$\{\varphi(q): q\in S(p)\}$ is a $\mu$-uniform cover of
$\varphi(p)$, i.\ e., there is a cover $\c{U}\in \mu$
such that $\c{U}\upharpoonright\varphi(p)\prec\varphi(S(p))$, or
in still other words, if
$ \varphi(S(p))\in \mu\upharpoonright\varphi(p)$. %%!! [S(p)] --> (S(p))
The principle introduced by the third author in
\referto[Pelanta] states that $\c{U}\in \lambda\mu$
if, and only if, there is a Noetherian tree
$T$ and a mapping $\varphi: T\to 2^X$, satisfying the
uniform mapping condition with respect to $\mu$, such
that $\varphi(\End(T))\prec\c{U}$. We may
delimit our discussion to trees consisting of
subsets of $X$, and define the uniform refinement
condition with respect to the associated inclusion
of the tree in $2^X$.

{\bf Perverse products.} We will define special `weak' subproducts
of direct products of trees with the property that
(under suitable conditions) infinite chains project onto
infinite chains in the factor trees.
In this paper, a {\it perversity}
is a decreasing sequence $p:\Bbb N\to \Bbb Z$ which is
eventually zero. (This concept is modified from
that used in intersection homology theory, cf.\ \referto[Kirwana].)
Sets of perversities, ordered by the coordinatewise
partial order of functions $\Bbb N\to \Bbb Z$,
are used to reduce a given direct product
of countably many trees.
Let $(T_i: i\in \Bbb N)$ be a countable family of trees,
and let $\c{P}$ be a set of perversities.

The
perverse product $\Pi_{\c{P}} T_i$ is the subset of the direct
product $\Pi_{i\in I} T_i$ consisting of all elements $x$ for which
there is a perversity $p\in \c{P}$ such that for each
$i\in \Bbb N$, the {\it level} of $x_i$ in $T_i$,
i.\ e., the number of predecessors in the tree, equals $p(i)$.
In addition to countable products, we define the finite
perverse products in the same way, by considering the
restriction $\c{P}\upharpoonright F = \{p\upharpoonright F: p\in \c{P}\}$
of perversities to the given finite set $F\subset \Bbb N$.
The {\it set-theoretic direct product} of the $T_i$ has as its
elements the Cartesian products $\Pi_i p_i$ (where each
factor $p_i$ is considered a set), and the
set-theoretic perverse product is the corresponding
subproduct.

The perverse product of countably many trees is a tree,
provided that the associated family of perversities is itself
a tree under its natural partial order.
The perverse product of
trees is a partially ordered set such that
infinite chains project onto
infinite chains in the factors, {\it provided that
the set of perversities has the same property}.

We now define a standard set $\c{P} = \{p_n: n\in \Bbb N\}$
of perversities to be used
in the sequel. Let $p_1 = (0,0,0,\ldots)$.
Suppose that
$p_n = (i_1,\ldots,i_m,0,0,\ldots)$ has been defined. Put
$$ k = \min\{{j\in\Bbb N}: i_j = i_{j+1}\}.$$
Then define $p_{n+1} = (i_1,\ldots,i_k+1,i_{k+1},\ldots,i_m,0,0,\ldots)$.
The elements of this sequence
are $()$, $(1)$, $(1,1)$, $(2,1)$, $(2,1,1)$, $(2,2,1)$ etc., where
we have indicated the non-zero entries only.
The set $\c{P}$ is linearly ordered with respect to the
coordinatewise order, and hence any perverse
product of trees with respect to $\c{P}$ is again a tree.
To give an example of its use,
let $T$ be a tree, and let
$T' = (\Pi_{i\in\Bbb N}T)_{\c{P}}$ be the perverse set-theoretic power
of $T$ with respect to $\c{P}$ as explained above.
The root of $T'$ (the unique element of level 0)
is the Cartesian power $R^{\Bbb N}$ of the root of $T$,
and the elements of level 3, say,
are products
$$ P_1\times P_2\times R\times R\times\ldots,$$
where $P_1$ is of level 2 and $P_2$ is of level 1 in $T$.

{\bf $\c{K}$-scattered spaces and $\c{K}$-derivatives.}
Let $\KK$ be a class of topological spaces. A
space $X$ is called {\it $\KK$-scattered}
\referto[Telgarskya],\referto[Telgarskyb] if every non-empty
closed subset of $X$ contains a point having a $\KK$-neighbourhood
in this subset.
We define the {\it $\KK$-derivatives} of a
$\KK$-scattered space $X$ as follows.
Let $D_{\KK}(X)$ denote the subset of all $p\in X$
having no $\KK$-neighbourhood. Then $D_{\KK}(X)$ is a closed
subset of $X$, and we set
$ D^{(0)}_{\KK}(X)= X$. Inductively, set
$D^{(\alpha+1)}_{\KK}(X)= D_{\KK}(D^{(\alpha)}_{\KK}(X))$,
and for a limit ordinal $\beta$ set
$D^{(\beta)}_{\KK}(X) = \cap_{\alpha<\beta} D^{(\alpha)}_{\KK}(X)$.
The $\c{K}$-rank of $X$ is defined as the smallest ordinal
$\alpha$ for which $D^{\alpha}_{\KK}(X) = \emptyset$.
Thus, the $\KK$-rank of a {\it locally} $\KK$ space is
at most one.

{\bf Decomposition trees.}
We will assume that $\KK$ is
hereditary with respect to closed subsets. The
{\it $\KK$-decomposition tree} $T_{\KK}(X)$ of $X$
will be defined as follows.
The elements of $T_{\KK}(X)$ are closed subsets of $X$
defined by using the derivative sets $D^{(\alpha)}_{\KK}(X)$.
Let $\bar\alpha = \Rank_{\KK}(X)$. Thus,
$D^{\bar\alpha}_{\KK}(X) = \emptyset$,
but $D^{\alpha}_{\KK}(X)\neq\emptyset$
for all $\alpha<\bar\alpha$. The closed subset
$\bigcap\{D^{\alpha}_{\KK}(X): \alpha<\bar\alpha\}$ is denoted by
the symbol $\top_{\KK}(X)$; notice that this is a locally
$\KK$ subset such that
every point $p$ of $X - \top_{\KK}(X)$ has a
closed neighbourhood $\bar U_p$ such that $\Rank_{\KK}(\bar U_p)
< \bar\alpha$. Also note that $\bar\alpha$ is a limit ordinal if,
and only if, $\top_{\KK}(X)$ is empty.
We let $\Root(T_{\KK}(X)) = X$, and the set
of immediate successors of $X$ is defined to consist of the set
$\top_{\KK}(X)$ and of {\it all} such closed subsets
$\bar U$ of $X - \top_{\KK}(X)$ with non-empty interior
and for which $\Rank_{\KK}(\bar U)
< \bar\alpha$.
If $\Rank_{\KK}(X) = 1$, then $T_{\KK}(X) = \{X\}$;
otherwise, the sets $\bar U$ satisfy $\Rank_{\KK}(\bar U)
<\bar\alpha$, and we can (recursively) define
the tree $T_{\KK}(X)$ by hanging
$\top_{\KK}(X)$ and
the trees $T_{\KK}(\bar U)$ below $X$.
Then $T_{\KK}(X)$ is a well-founded tree.

\section{Partition-completeness and stationary winning strategies.}
\markplace
%%!!E.\ Michael defined in \referto[Michaelb] a generalization of
%%!!Frol\'\i k's and Arhangel'skii's complete sequences
%%!!\referto[Frolikb], \referto[Arhangelskia] that
%%!!characterize \v Cech-complete spaces. 
In this section, we consider a condition that 
generalizes both C-scattered and \v Cech-complete spaces.
Recall that a
sequence $\hbox{\bs U}=(\c{U}_n)$ of
covers of a space $X$ is {\it complete}
if any filter base $\c{F}$ {\it `controlled'} by $\hbox{\bs U}$
has a cluster point, in other words if any filter $\c{F}$
with $\c{F}\cap \c{U}_n\neq\emptyset$ for all $n$ satisfies
$\bigcap\{\bar F: F\in\c{F}\}\neq\emptyset$. 
The classical characterization from Frol\'\i k\referto[Frolikb] and 
Arhangel'skii\referto[Arhangelskia] describes the \v Cech-complete 
spaces as the completely regular spaces with a complete sequence 
of open covers. 
%%!!A cover $\c{U}$ of $X$
%%!!is called {\it exhaustive} if every non-empty subset $S$ of $X$
%%!!has a non-empty relatively open subset of the form $U\cap S$,
%%!!where $U\in \c{U}$. Michael's condition generalizing
%%!!that of Frol\'\i k and Arhangel'skii is that the
%%!!space have a complete sequence of exhaustive covers.

On the other hand, $\c{K}$-scattered spaces have canonical
exhaustions into {\it left-open partitions} of subsets
$S$ such that $\bar S$ is a $\c{K}$-subset. 
%%!!(For example,
%%!!\referto[Hohtid], \referto[HohtiPelantb] considered
%%!!exhaustions of C-scattered spaces rather than their
%%!!inverted Noetherian decomposition trees.)
A partition $\c{P}$ of
a space $X$ is called {\it left-open} if $\c{P}$ admits a
well-ordering $\{P_\alpha: \alpha < \beta\}$ such that
$\bigcup\{P_\alpha: \alpha \leq \gamma\}$ is open for all
$\gamma < \beta$. It was proved in \referto[JunnilaSmithTelgarsky]
that if $\c{K}$ is closed-hereditary, then $X$ is
$\c{K}$-scattered if and only if $X$ has a left-open
partition $\c{P}$ such that $\bar P$ is a $\c{K}$-subset
for all $P\in \c{P}$.
In \referto[TelgarskyWicke],
Telg\'arsky and Wicke studied spaces which have a
complete sequence of left-open partitions: A sequence
$\hbox{\bs P}=(\c{P}_n)$ of left-open partitions is
called {\it complete} if 1) $\c{P}_{n+1}$ refines
$\c{P}_n$ for all $n$, 2) the well-order of $\c{P}_{n+1}$
is compatible with that of $\c{P}_n$ in the sense that
the relation $Q\subset P$ between elements
$P\in \c{P}_n$, $Q\in \c{P}_{n+1}$ preserves the order, and
3) any filter base controlled by $\hbox{\bs P}$ has
a cluster point. 

\th Definition {\ThisSection}.1 (\referto[TelgarskyWicke], p.\ 737): A space 
is called partition-complete if it has a complete 
sequence of left-open partitions.

We call a space $\sigma$-partition-complete if it is the
union of countably many partition-complete, closed subspaces.

The results important for us were already obtained 
by E.\ Michael \referto[Michaelb].
Recall that a space is $\c{K}$-scattered if every non-empty
closed subset of the space contains a point with a $\c{K}$-neighbourhood
in the subset. 
%%!!This is equivalent to the following: Every non-empty
%%!!subset $S$ contains a relatively open subset $U$ of the form
%%!!$K\cap S$, where $\bar K$ is in $\c{K}$ ($\c{K}$ is closed-hereditary).
Michael called a cover $\c{U}$ of $X$
{\it exhaustive} if every non-empty subset $S$ of $X$
has a non-empty relatively open subset of the form $U\cap S$,
where $U\in \c{U}$. (Thus, for C-scattered spaces, the 
interiors of compact subsets form an exhaustive cover.)
The following result of Michael gives us the definition of 
partition-completeness we will use in this paper. (The result
follows directly from \referto[Michaelb], Prop.\ 4.1 and 
\referto[TelgarskyWicke], Prop.\ 1.2.)

\th Theorem {\ThisSection}.2 (\referto[Michaelb], Prop.\ 4.1): A space
$X$ is partition-complete if, and only if, it has a complete sequence of 
exhaustive covers.

%%!!Such {\it partition-complete} spaces are
%%!!by \referto[Michaelb], Prop.\ 4.1, the same as the spaces
%%!!$X$ with a complete exhaustive cover, and we will
%%!!consider them with respect to the stationary winning
%%!!strategies $\phi$ for Player II in the game $G(X)$ defined above.

We will use this condition in terms of a topological
game $G(X)$ related to such sequences of covers.
Players I and II choose in alternating steps non-empty
subsets $S_1\supset T_1\supset S_2\supset\cdots$
of $X$ such that Player I chooses the sets $S_i$,
Player II chooses the sets $T_i$, and $T_i$ is
relatively open in $S_i$ for all $i$. Player II wins
if any filter base finer than $\{T_n: n\in \Bbb N\}$
has a cluster point. A {\it stationary winning strategy}
for Player II in $G(X)$ is simply a function
$\phi: P(X)\to P(X)$ giving the relatively open choices $\phi(S)\subset S$ 
($\phi(S)\neq\emptyset$)
of Player II for all
$S\neq\emptyset$.
Michael proved that a space $X$
has a complete sequence of exhaustive covers if and only
if Player II has such a stationary winning strategy in
the game $G(X)$ (\referto[Michaelb], Prop.\ 4.1 and
Thm.\ 7.3).

Notice the connection between $\c{K}$-scattered spaces and
stationary winning strategies in the above sense. By
providing for each non-empty closed subset $F\subset X$
an open set $\emptyset\neq\phi(F)\subset F$, the
winning strategy $\phi$ yields a decreasing decomposition
of $X$ into `derivative' subsets and thus (inverting the order)
into a Noetherian
decomposition tree. The $\phi$-derivative of a non-empty subset
$S\subset X$ is simply defined as the relatively closed set $S\setminus \phi(S)$.
The consecutive derivatives are then defined in a complete
analogy with the definition given in the previous section.
We denote the corresponding decomposition tree by $T_\phi(X)$.

If $X$ is partition-complete with respect to a winning strategy
$\phi$, we extend its decomposition tree $T_\phi(X)$ to
a (in general) non-Noetherian `game tree' as follows. For each
$U\in \End(T_\phi(X))$, let $\c{A}(U)$ be the collection of
all closed $F\subset U$ of the form $\bar V$, where
$V\subset U$ is non-empty and open. In order to ensure that
$U$ is the union of such closures, we assume that
$X$ is regular. Denote by $\phi_S$ the restriction of
$\phi$ to the subsets of $S$, for any
subset $S\subset X$. Extend $T_1 = T_\phi(X)$ to
$T_2$ by declaring the $\bar V\in \c{A}(U)$ as the immediate
successors of $U$, and hanging the trees $T_{\phi_{\bar V}}(\bar V)$
below $\bar V$. This operation is repeated countably
many times; we obtain an increasing sequence
$T_1,T_2,\ldots$ of Noetherian trees, and define
$T(X)$ as their union. Although $T(X)$ in general is
not well-founded, it has the following Noetherian
property with respect to open covers. If $\c{G}$ is an
open cover of $X$, then there is a Noetherian subtree
$T'$ of $T(X)$ such that $\End(T') \prec \c{G}^{<\omega}$,
where $\c{G}^{<\omega}$ denotes the directed open cover
obtained as the finite unions of elements in $\c{G}$.
This follows from the following general (and simple) lemma:

\th Lemma {\ThisSection}.3: Let $X$ be a topological
space, and let $ \c{F} = (F_n: n\in \Bbb N)$ be a decreasing sequence of
closed subsets such that any filter finer %%!! filter base
than $\c{F}$ has a cluster point. Then $\bigcap\{F_n: n\in\Bbb N\}$ is
compact.

\section{The $\lambda$-neighbourhood induction lemma.} Our main method
for proving that a given cover $\c{G}$ belongs to the locally fine
coreflection applies the following `Noetherian' induction: There is
a Noetherian tree $T$ such that each maximal element $E$ of $T$
has a $\lambda$-neighbourhood $U$ with the property that
$\c{G}\upharpoonright U$ is a $\lambda$-cover of $U$. Then the
construction of $T$ as a (subtree of) perverse power and the
induction lemma given below enable us to conclude the
validity of this property for the top element of $T$. This
was stated and proved in \referto[HohtiYun] for uniform spaces,
but it is valid for pre-uniformities (covering monoids).

\th Lemma {\ThisSection}.1: Let $X$ be a set, and let
$\mu$ be a pre-uniformity on $X$. Let $A\subset X$, let
$\c{V}\in \lambda\mu\upharpoonright A$ and for each $V\in \c{V}$,
let $V'$ be a $\lambda\mu$-uniform neighbourhood
of $V$ in $X$ with respect to $\mu$. Then $\bigcup\{V': V\in \c{V}\}$
is a $\lambda\mu$-uniform neighbourhood of $A$ in $X$. Moreover,
if $\c{G}$ is a $\lambda\mu$-uniform cover of each $V'$, then
$A$ has a $\lambda\mu$-uniform neighbourhood $N$ such that
$\c{G}$ is a $\lambda\mu$-uniform cover of $N$.

\section{The Metric-fine Coreflection.} We call a cover $\c{G}$
of a uniform space $\mu X$ {\it $\sigma$-uniform} if there is
a countable closed cover $\{F_n: n\in \Bbb N\}$ of $X$ such that
each restriction $\c{G}\upharpoonright F_n$ is a uniform cover
of the subspace $F_n$. This property can be formulated in
several other settings. One calls \referto[Hagera] a uniform
space $\mu X$ {\it metric-fine}, if for any metric uniform space
and any uniformly continuous map $f: \mu X\to \rho M$, $f$
remains uniformly continuous when the target $\rho M$ is
changed to $\c{F}M$, i.\ e., the fine space associated with
$M$. The {\it metric-fine coreflection} $m\mu$ is the weakest
metric-fine uniformity stronger than $\mu$. We note that $m\mu$ is
is the collection of all $\sigma$-uniform
covers of $X$. Each $\sigma$-uniform (open) cover has a
$\sigma$-uniformly discrete (open) refinement. A cover $\c{G}$
is $\sigma$-uniform if, and only if, there is a sequence
$(\c{U}_n: n\in \Bbb N)$ of uniform covers of $\mu X$ such that
for each $x\in X$, there is $n(x)\in \Bbb N$ with
$\st(x,\c{U}_{n(x)})\prec \c{G}$. (Notice the relation with the
so-called {\it $\theta$-refinability} of (topological) spaces).
The metric-fine coreflection was used in both \referto[HohtiPelantb]
and \referto[HohtiYun] to extend the results from $\c{K}$-scattered
spaces to $\sigma-\c{K}$-scattered, where $\c{K}$ was the
class of compact (resp.\ \v Cech-scattered paracompact) spaces.
The metric-fine coreflection $m$ can be
directly extended to pre-uniformities of topological spaces.
In our situation, we define the metric-fine modification 
$m\mu$ of a pre-uniformity $\mu$ as the filter of 
all covers refined by covers $\c{V}$ for which there is 
a countable closed cover $\{F_n: n\in \Bbb N\}$ of the 
underlying space such that for each $n\in \Bbb N$, there 
is $\c{U}_n\in \mu$ with 
$\c{U}_n\upharpoonright F_n \prec \c{V}\upharpoonright F_n$.

The covers in $\lambda m\mu$ are obtained by means of Noetherian
trees as follows. There is a Noetherian tree $T$ and mapping
$\varphi: T\to 2^X$ with the
property that for each $p\in T$, the immediate
successors $S_p$ of $p$ in $T$ are divided into countably
many parts $S^{(n)}_p$ such that for each $n$, the image
$\varphi(S^{(n)}_p)$ is a uniform cover of the closed subspace
$\bigcup(\varphi(S^{(n)}_p))$ of $\varphi(p)$, and these
subspaces form a countable cover of $\varphi(p)$.

%%!! This section moved to the end ("Additional Remarks")
%%!!\section{Frames, formal topology and `effective' refinements.}

\section{Countable Products of $\sigma$-partition-complete
supercomplete spaces.}

\markplace{FirstMainSection}{\ThisSection}

In our proof of \refered[ThNormalCovers]  %%!!
that the locally fine coreflection of the
product $\lambda\Pi\c{F}X_i$ contains all the normal covers, we
will need the result that each finite
product is
{\ninepoint supercomplete\footnote{\dag}{Nevertheless, it is not
sufficient that all the finite products are supercomplete.
It was shown in \referto[Hohtie] that there is a subspace
$X$ of the real numbers such that each finite power
$(\c{F}X)^n$ is supercomplete but $(\c{F}X)^\omega$ is not.
The space $X$ is not partition-complete. On the other
hand, it was shown in \referto[Hohtif] that there are
non-analytic subspaces $X$ of $\Bbb R$ such that
$\lambda((\c{F}X)^\kappa) = (\c{F}X)^\kappa$ for
all $\kappa$.}}.
We establish here a stronger result for countable products.
Indeed, we don't have to stop at paracompactness. We will
prove a product theorem for regular spaces $X$ having
{\it supercomplete pre-uniformities} (monoids of covers),
i.\ e., filters $\mu$ of coverings such that each open
cover of $X$ is in $\mu$. We first state the result
for partition-complete spaces.

\setbox0=\hbox{\number\sectioncount.1\ }
\def\CPTHM{\copy0}
\th Theorem {\ThisSection}.1: Let $(X_i: i\in \Bbb N)$
be a countable family of partition-complete, regular
spaces, and for each $i$, let $\mu_i$ be a supercomplete pre-uniformity
for $X_i$. Then
$\lambda\Pi\mu_i X_i$ is a supercomplete pre-uniformity for
the topological product $\Pi X_i$.

\proof The proof is very similar to the one given in
\referto[HohtiYun] for the case of countably many
\v Cech-scattered paracompacta. Therefore, we will only
give the main steps of the proof and refer the reader to
\referto[HohtiYun] for details. Let $\c{G}$ be an open cover
of $\Pi X_i$.
We may assume that the
factors $X_i$ are the same space $X$. (Replace the factors
by their disjoint union.)
For each $i$, let
$T$ be the tree $T_\phi(X)$, where
$\phi$ is the stationary winning strategy. Let
$T'$ be the perverse power $\Pi_{\c{P}}T$, where
$\c{P}$ is the standard set of perversities defined above (Section 2). %%!!
Finally,
let $T''$ be the subtree of $T'$ consisting of all
$P\in T'$ which do not have a predecessor $Q$
with a uniform neighbourhood $N$ such that $\c{G}$ is
a uniform cover of $N$. (The term `uniform' refers
here to the related pre-uniformities.)

We claim that $T''$ is well-founded. Indeed, suppose to the
contrary that $T''$ has an infinite chain $(P_n: n\in \Bbb N)$.
By the definition of perverse products, the projections
$(\pi_i[P_n]: n\in \Bbb N)$ are again sequences with
infinitely many distinct terms.
Thus there are, for each $i$, relatively open sets
$E^{(i)}_n\subset X$ such that
$$ \ldots\supset E^{(i)}_n \supset \pi_i[P_n] \supset
               E^{(i)}_{n+1} \supset \ldots,$$
indeed, we may write $E^{(i)}_{n+1} = \phi(\pi_i[P_n])$.
Therefore, the intersection
$$ K_i = \bigcap_{n\in \Bbb N} \pi_i[P_n]$$
is compact, and it follows that
$$ K = \bigcap_{n\in \Bbb N} P_n = \Pi_{i\in \Bbb N} K_i$$
is compact, too.

For each $x\in K$, choose a basic open set
$B_x$ such that $x\in B_x\subset \bar B_x\subset G_x$
for some $G_x\in \c{G}$. By the compactness of $K$,
there is a finite set $\c{B}=\{B_x : x\in F\}$ which covers
$K$. It is easy to find (finite) families $\c{B}_i$ of open
subsets of the spaces $X_i$ which cover the sets $K_i$
such that for some $s\in \Bbb N$,
$$ \c{B}' \;=\; \bigcap_{i=1}^s\pi_i^{-1}[\c{B}_i] \prec \c{B}.$$
Furthermore, by the regularity of $X_i$
we can assume there are open neighbourhoods
$N_i, M_i$ of $K_i$ covered by $\c{B}_i$ such
that $\bar N_i\subset M_i$.
As every open cover of
$X_i$ is in $\mu_i$, it follows that $\c{B}'$ is a uniform
cover of the uniform neighbourhood
$$ N \;=\; \bigcap_{i=1}^s \pi_i^{-1}[N_i]$$
of $K$, and therefore so is $\c{G}$. Now for each
$i\in \{1,\ldots,s\}$ there is $n_i$ such that
$\pi_i[P_{n_i}]\subset N_i$. Put
$n = \max\{n_1,\ldots,n_s\}$. Then it is easy to see that
$ P_n\subset N$,
which yields a contradiction, because $P_{n+1}$ is an
element of $T''$.

We prove by using the $\lambda$-neighbourhood induction lemma
that $\c{G}$ is a $\lambda$-uniform cover of the root of $T''$. Let
$P$ be an element of $T''$. If $P$ is a maximal element
of $T''$, then by the definition of $T''$, there is a
uniform neighbourhood $N_P$ of $P$ such that $\c{G}$ is a
uniform cover of $N_P$. (More precisely, there is a
cover $\c{U}$ in $\Pi\mu_i$ such that
$\c{U}\upharpoonright N_P \prec \c{G}\upharpoonright N_P$.)
We proceed by induction and assume that every successor
$Q$ of $P$ in $T''$ has a $\lambda$-uniform neighbourhood
$N_Q$ such that $\c{G}$ is a $\lambda$-uniform cover of
$N_Q$. Since $T'$ is a perverse product, and since
$T''$ is a lower part (an ideal) of $T'$, there is
$i\in \Bbb N$ such that the immediate successors
$Q$ of $P$ differ only with respect to $i$: For
$j\neq i$, $\pi_j[P] = \pi_j[Q]$, while $\pi_i[Q]$ is an immediate
successor of $\pi_i[P]$ in $T_\phi(X_i)$. We recall that
$T(X)$ was obtained through a countable iteration of the
trees $T_\phi(S)$ for closed subsets $S$ of $X$, defining it
as the union of the trees $T_n$. There are two
possibilities. Either $\pi_i[P]$ is a maximal element
of some $T_n$ and its immediate successors all belong
to $T_{n+1}$, or else both $\pi_i[P]$  and its immediate
successors belong to the same tree of the form $T_\phi(S)$.
However, in each case, the immediate successors are divided
into two parts: a closed set $F\subset \pi_i[P]$ and
all the closures $\bar V$ of non-empty open sets
$V\subset \pi_i[P]\setminus F$
such that $\bar V\cap F = \emptyset$.

Let $Q_F = \pi_i^{-1}[F]\cap P$ be the immediate successor
of $P$ that corresponds to $F$. By the induction hypothesis,
$Q_F$ has a $\lambda$-uniform neighbourhood $N_F$ such that
$\c{G}$ is a $\lambda$-uniform cover of $N_F$.
Then $P - N_F$ has a closed $\lambda$-uniform
neighbourhood $V_F$ such that
$V_F\cap Q_F = \emptyset$, and such that the binary cover
$\{N_F,V_F\}$ is a $\lambda$-uniform cover of $P$.
Indeed, there is a $\lambda$-uniform cover
$\c{V}$ of $X^{\omega}$ such that
$\st(P-N_F,\c{V})\cap Q_F = \emptyset$.
(Simply take $\c{V}$ as a $\lambda$-uniform cover such that
$\st(Q_F,\c{V})\subset N_F$.)
We can assume that $\c{V}$
consists of the closures of basic open sets.
This follows from the assumption that the spaces
$X_i$ are regular and that every open cover of $X_i$
belongs to $\mu_i$.
(Recall that the locally fine coreflection
of a product of pre-uniformities has a basis
of covers consisting of rectangular elements.)
Then each element $V$ of
$\c{V}\upharpoonright P$
that meets $P-N_F$ satisfies $\pi_i[V]\cap \pi_i[Q_F] = \emptyset$,
which implies that $V$ is an immediate successor of $P$. But then
by the inductive hypothesis, each $V\in\c{V}\upharpoonright P$ with
$V\cap (P-N_F)\neq\emptyset$ has a $\lambda$-uniform neighbourhood
$N_V$ such that $\c{G}$ is $\lambda$-uniform cover of $N_V$.
It then follows that $P$ has a $\lambda$-uniform neighbourhood
$N_P$ satisfying the condition of our claim.
In fact, by the $\lambda$-neighbourhood induction lemma
the set $N_F$ has a $\lambda$-uniform
neighbourhood $N_F'$ such that ${\cal G}$
is a $\lambda$-uniform cover of $N_F'$.
The sets $V$ such that
$V\cap (P\setminus N_F)\neq\emptyset$ together with
$N_F$
form a $\lambda$-uniform cover of $P$. As ${\cal G}$
is $\lambda$-uniform over the sets $N_V$ and
$N_F'$, the set $P$ has by a new application of
the lemma a $\lambda$-uniform neighbourhood $N_P$ with the
desired property. $\eop$

For the following corollaries, notice that the pre-uniformities
generated by all %%!! formed by
locally finite, locally countable, point-finite, point-countable,
and countable open covers are locally fine (i.\ e., stable
under the operation $\lambda$).

\th Corollary {\ThisSection}.2: Let $(X_i: i\in \Bbb N)$
be a countable family of regular paracompact (resp.\ para-Lindel\"of,
metacompact, meta-Lindel\"of, Lindel\"of) spaces. If the spaces
$X_i$ are partition-complete, then the product space
$\Pi X_i$ is paracompact (resp.\ para-Lindel\"of, metacompact,
meta-Lindel\"of, Lindel\"of).

\proof Suppose that the spaces $X_i$ are paracompact and partition
complete. Let $\mu_i$ be the pre-uniformity formed by all locally
finite open covers of $X_i$. Then $\mu_i$ is supercomplete for
$\Omega(X_i)$, because by paracompactness every open cover
of $X_i$ can be refined by, and thus is, a member of $\mu_i$.
By the theorem proved above, $\lambda\Pi\mu_i$ is 
supercomplete  pre-uniformity for the topology of $\Pi X_i$.
But $\Pi\mu_i$ has a basis consisting of products of finitely many
locally finite open covers, and therefore so does
$\lambda\Pi\mu_i$, because the 
pre-uniformity generated by the locally finite open covers %%!! 
is preserved by $\lambda$. It follows that every open cover %%!! are
of $\Pi X_i$ has an open locally finite refinement, and thus
the product is paracompact. The other cases are proved
in an analogous way. $\eop$

This corollary can be stated also for other locally fine
covering concepts, e.\ g., ultraparacompactness (refinement
by clopen covers) etc. It was essentially proved by
Plewe in \referto[Pleweb]. He showed that
the localic product of countably many partition-complete
regular spaces is spatial. The result then follows
from a theorem of Dowker and Strauss (\referto[DowkerStrauss])
that the localic product of regular paracompact (resp.\
metacompact, Lindel\"of, ultraparacompact) spaces
is again paracompact (resp.\ metacompact, Lindel\"of,
ultraparacompact). The case of paracompactness was already
due to Isbell \referto[Isbellc]. One could also
obtain theorem \CPTHM by combining Plewe's
result on spatiality and the relationship
between $\lambda\Pi\mu_i$ and $\otimes\Omega(X_i)$,
mentioned above. In addition to the
topological corollaries,
theorem \CPTHM yields stronger combinatorial
corollaries\footnote{\dag}{\ninepoint We also note that
it gives a basis for effective product theorems, inasmuch
as the factor spaces can be presented effectively.}
for uniformities.

\th Corollary {\ThisSection}.3: Let $(\mu_iX_i: i\in \Bbb N)$
be a countable family of supercomplete
spaces. If the spaces
$X_i$ are partition-complete, then the product space
$\Pi\mu_i X_i$ is supercomplete.

\proof As the spaces $\mu_iX_i$ are supercomplete, we have 
$\lambda\mu_i = \c{F}(X_i)$ for each $i$. Therefore, by
{\ThisSection}.1, 
 $$\lambda \Pi\mu_i = \lambda \Pi \lambda \mu_i =
 \lambda \Pi {\cal F}(X_i) = \c{F}(\Pi X_i),$$
where we have used \referto[Isbella], Exercise VII 8a (p.\ 143).

As in \referto[HohtiPelantb], \referto[HohtiYun], we may extend 
the proof of \CPTHM to the $\sigma$-partition-complete case.
(Recall that a space is called $\sigma$-partition-complete if
it is a countable union of partition-complete, closed subspaces.) %%!!
Thus, we obtain the following result: %%!! move proof after Theorem

\markplace{ThSigmaPartitionComplete}{Theorem\ \ThisSection.4}

\th Theorem {\ThisSection}.4: Let $(X_i: i\in \Bbb N)$
be a countable family of $\sigma$-partition-complete, regular
spaces, and let $\mu_i$ be a supercomplete pre-uniformity
for $X_i$. Then
$\lambda m\Pi\mu_i$ is a supercomplete pre-uniformity for
the topological product $\Pi X_i$.

\proof Given a family $(X_i: i\in \Bbb N)$ of regular spaces with 
supercomplete pre-uniformities $\mu_i$,
each of which is a countable union of closed, partition-complete
subspaces, we replace the trees $T_\phi(X_i)$ with a forest which 
contains a subtree $T_\phi(F_n)$ for each of the partition 
complete parts $F_n$. 
The perverse product construction 
ensures that only products of partition-complete parts are used.
Since the partially ordered set $T''$ so constructed is Noetherian, $T''$ only has 
countably many 'types' of subproducts.
We use the same induction as 
in \CPTHM, but we arrive, instead of the 
entire product of the $X_i$, at minimal elements $P$
which are subproducts and have ($\lambda$-uniform) 
neighbourhoods $N_P$ such that $\c{G}^{\omega}$ is 
$\lambda$-uniform over $N_P$. As the successors
of these minimal elements belong to the same type
(i.\ e., their factors are contained in the same 
partition-complete parts), and there are only 
countably many types, we arrive at a countable, closed cover
of the product such that $\c{G}$ is $\lambda$-uniform over 
each member. Consequently, $\lambda(m\Pi\mu_i)$ is a 
supercomplete pre-uniformity. 

\section{Normal covers in uncountable products.} In \referto[Pelanta],
the third author proved that every normal cover of an
arbitrary product of complete metric spaces belongs to the
locally fine coreflection of the corresponding product of
the fine uniformities of the factor spaces. This result had
been previously established for arbitrary products of
Polish spaces (cf.\ \referto[Isbellb], VII, Cor.\ 21.). In this
section we will extend this to products of paracompact,
partition-complete spaces. It turns out as in the preceding
section that we may consider supercomplete pre-uniformities
of regular partition-complete spaces.

The proof is based on simple combinatorial properties of
`basic subsets' of products. Let $(X_i: i\in I)$  be a
family of topological spaces. The basic sets $B\subset \Pi X_i$
are of the form $B = \cap\{\pi_i^{-1}[B_i]: i\in F\}$, where
$B_i\subset X_i$ and $F\subset I$ is finite. The
{\it basic open sets} are those for which each $B_i$ is open
in $X_i$. Given a subset $S\subset \Pi X_i$,
denote by $\c{B}(S)$ the collection of all basic subsets $B$
of $\Pi X_i$ such that $B\subset S$.
If $B\in \c{B}(\Pi X_i)$, then $I(B)$ is the set of
all $i\in I$ such that $\pi_i[B]\neq X_i$.

We first recall lemmas from \referto[Pelanta].
It is
useful to begin with the following:

\th Lemma {\ThisSection}.1:
Let $B_1$, $B_2$ be basic
sets in a product $\Pi X_i$, and let
$F = I(B_1)\cap I(B_2)$. Suppose that $B_1,B_2\neq\emptyset$.
\item{1)} If $F=\emptyset$, then $B_1\cap B_2 \neq\emptyset$.
\item{2)} If $F\neq\emptyset$, then $B_1\cap B_2\neq\emptyset$ if
and only if $\pi_F[B_1]\cap \pi_F[B_2]\neq\emptyset$.  %%!!

\th Lemma {\ThisSection}.2 (\referto[Pelanta]):
If $\{B_n: \in \Bbb N\}$ is a
family of non-empty basic sets in $\Pi X_i$ such that the index
sets $I(B_n)$ are mutually disjoint, then
$\overline{\cup\{B_n: n\in \Bbb N\}} = X$.

\proof Indeed, for a point $x\in X$ and for a basic open
neighbourhood $B$ of $x$ there is $n$ with $I(B)\cap I(B_n)=\emptyset$.
By the preceding lemma, this implies $B\cap B_n\neq\emptyset$,
and thus $x$ belongs to the closure of the union of the $B_n$. $\eop$

The following is the essential {\it inclusion lemma} for regularly
open subsets of products:

\th Lemma {\ThisSection}.3 (\referto[Pelanta]): Let $G$
be an open subset and let $R$ be a regular open subset of $X=\Pi_I X_i$.
Suppose that there are $E\subset I$ and a sequence $(B_n: n\in \Bbb N)$
of basic subsets of $R$ such that $\pi_E[B_n]\supset \pi_E[G]$ for
all $n$, and that the index sets $I(B_n)\setminus E$ are mutually
disjoint. Then $G\subset R$.

\th Lemma {\ThisSection}.4: Let $R$ be a regular open
proper subset of $X = \Pi_I X_i$. Then there is a finite set
$F(R)\subset I$ which meets $I(B)$ for each $B\in \c{B}(R)$.

\proof Suppose that the condition of the lemma is not
satisfied. Thus, for each finite subset $\c{B}'\subset \c{B}(R)$
there is $B'\in \c{B}(R)$ such that
$I(B')\cap \cup\{I(B): B\in \c{B}'\}=\emptyset$. Therefore,
$\{I(B): B\in \c{B}(R)\}$ contains an infinite
disjoint subfamily, and by Lemma {\ThisSection}.2
we have $R = X$, which contradicts the assumption. $\eop$

\def\extended-lemmas
{
%%%%%%%%%%%%%%%%%%%%%%%%%%%%%%%%%%%%%%%%%%%%%%%%%%%%%%%%%%%%%%
%%%%%%%%%%%%%%%%%%%%%%%%%%%%%%%%%%%%%%%%%%%%%%%%%%%%%%%%%%%%%%
In this paper, we will use two extended forms of these lemmas
for finite unions of basic sets. They follow from the above
versions.

%\th Lemma {\ThisSection}.1':
%Let $\{B_i: i\in F_1\}$, $\{B_i: i\in F_2\}$ be finite
%families of basic
%sets in $X =\Pi_I X_i$, and let
%$ G_1 = \cup\{B_i:i\in F_1\}$ and $G_2 = \cup\{B_i:i\in F_2\}$.
%If there is $E\subset I$ such that
%$\pi_E[G_1]\cap \pi_E[G_2]\neq\emptyset$ and
%$$\left(\cup_{i\in F_1} I(B_i)\right)\cap
%                  \left(\cup_{i\in F_2} I(B_i)\right)\subset E,$$
%then $G_1\cap G_2\neq\emptyset$.
%
%\proof It follows from the first assumption that
%there are $i\in F_1$, $j\in F_2$ with
%$\pi_E[B_i]\cap \pi_E[B_j]\neq\emptyset$. In addition,
%$I(B_i)\cap I(B_j)\subset E$. By Lemma
%{\ThisSection}.1, $B_i\cap B_j\neq\emptyset$.
%Thus, $G_1\cap G_2\neq\emptyset$. $\eop$

\th Lemma {\ThisSection}.1':
If $\{\c{B}_n: n\in \Bbb N\}$ is a
family of finite families of basic sets
in $X = \Pi_I X_i$ such that the index
sets $I(\c{B}_n) = \cup\{I(B): B\in \c{B}_n\}$
are mutually disjoint, then the union of the sets
$\cup(\c{B}_n)$ is dense in $X$.

\proof Indeed, this is weaker than the statement for
one-element collections. $\eop$

\th Lemma {\ThisSection}.2': Let $G$
be an open subset and let $R$ be a
regular open subset of $X=\Pi_I X_i$.
Suppose that there is a sequence $(\c{B}_n: n\in \Bbb N)$
of finite families of basic subsets of $R$
such that the index sets $I(\c{B}_n)$ are mutually disjoint
relative to a subset $E\subset I$ for which
$\pi_E[\cup(\c{B}_n)]\supset \pi_E[G]$ for all $n$. Then
$G\subset R$.

\proof By the assumption, and by the preceding lemma,
$$ G\subset \pi_E[G] \times \Pi_{I\setminus E} X_i \subset
           \overline{\cup_{n\in \Bbb N}\left(\cup(\c{B}_n)\right)},$$
because $\pi_E[G]\subset \cup\{\pi_E[B]: B\in \c{B}_n\}$ and the
basic sets $B$ satisfy the
condition $B = \pi_E[B]\times \pi_{I\setminus E}[B]$. Thus,
$G\subset \int(\bar R)$, and hence $G\subset R$, because
$R$ is regular open. $\eop$

From this we obtain the following corollary:

\th Lemma \ThisSection.3': Let $G,R$ be open subsets
of $X = \Pi_I X_i$, where $G$ is not contained in $\int(\bar R)$.
Let $E\subset I$ be finite. Then there is a finite set
$F(R)\subset I$ such that given any finite family $\c{B}$ of
basic open subsets $B\subset R$ with
$\pi_E[\cup(\c{B})]\supset \pi_E[G]$, we have
$$ (F(R)\cap (\cup\{I(B): B\in \c{B}\})) \setminus E \neq \emptyset. $$

\proof If no such $F(R)$ existed, we could inductively
construct a sequence $\{\c{B}_n: n\in \Bbb N\}$ of
finite families $\c{B}_n$ of basic open sets such that
$\pi_E[\cup(\c{B}_n)] \supset \pi_E[G]$ and for which
the index sets $I(\c{B}_n)$ are pairwise disjoint relative to $E$. But
then by the preceding lemma,
we would have $G\subset R^*$, contradicting the assumption. $\eop$

} %%%%%%%%%%%%%%%%%%%%%%%%%%%%%%%%%%%%%%%%%%%%%%%%%%%%%%%%%%
%%%%%%%%%%%%%%%%%%%%%%%%%%%%%%%%%%%%%%%%%%%%%%%%%%%%%%%%%%%%
\def\regular-open-extension{
For an open subset $S$ of a
topological space $X$, the {\it regular
open extension}
$S^*$ is the subset $\int(\bar S)$.
Suppose that $\c{V}$ is a normal
cover, i.\ e., $\c{V}$ is uniformizable. Let $\c{V}_1$ be a
normal cover with $\c{V}_1\prec^{**}\c{V}$, and let $\c{R}$ be a
locally finite, open refinement of $\c{V}_1$ (which
every normal cover has). Then the regular open extension
$\c{R}^* = \{R^*: R\in \c{R}\}$ of $\c{R}$ refines $\c{V}$.
Indeed, let $R\in \c{R}$, and let $R\subset V\in \c{V}_1$,
$\st(V,\c{V}_1)\subset V'\in \c{V}$. Then
$\int(\bar R)\subset \int(\bar V) \subset \st(V,\c{V}_1)\subset V'$.
Moreover, $\c{R}^*$ is locally finite.

In the following
proof, we will replace the given normal cover $\c{V}$
by the regular open extension of $\c{R}^{<\omega}$,
where $\c{R}^{<\omega}$ is the directed cover consisting
of the finite unions of elements of $\c{R}$, and $\c{R}$ is
as above. Indeed, let $\c{V}, \c{V}_1$ and $\c{R}$ be as above.
We claim that the regular open extension $(\c{R}^{<\omega})^*$ of
$\c{R}^{<\omega}$ refines $\c{V}^{<\omega}$. To see this,
let $R_1,\ldots, R_n\in \c{R}$, where we have sets $V_i'\in \c{V}_1$
such that $R_i\subset V_i$ and $\st(V_i',\c{V}_1)\subset V_i\in \c{V}$.
Then
$$ \int(\overline{R_1\cup\cdots\cup R_n}) \subset
      (V_1'\cup\cdots\cup V_n')^* \subset V_1\cup\cdots\cup V_n. $$
}

\extended-lemmas

\regular-open-extension

We are now able to state the main result of this section.

\markplace{ThNormalCovers}{Theorem\ \ThisSection.5}

\th Theorem {\ThisSection}.5: Let $(X_i: i\in I)$
be a family of partition-complete, completely regular spaces, and for
each $i\in I$, let $\mu_i$ be a supercomplete uniformity
for $X_i$. Then every normal cover $\c{V}$ of
$X = \Pi_I X_i$ belongs to $\lambda\Pi_I\mu_i$.

\proof
Let $\c{R}$ be a locally finite
cover of $X$ by regular open
sets such that $(\c{R}^{<\omega})^*\prec \c{V}^{<\omega}$,
and let $\c{G}$ be a cover consisting of basic open sets
such that each member of $\c{G}$ meets only finitely many
members from $\c{R}$. We will show that
$\c{V}^{<\omega}$ belongs to $\lambda\Pi_I\mu_i$.
We conclude that every normal cover belongs
to this product uniformity.

We may assume that
the cover $\c{R}$ is non-trivial, i.\ e.,
$R\neq X$ for all $R\in \c{R}$.
Given any $G\in \c{G}$, consider the
finite set $\c{R}_G\subset \c{R}$  such that $G\cap R = \emptyset$
for all $R\in \c{R}\setminus \c{R}_G$. Then $B\in \c{B}(R)$
for such an $R$ implies $G\cap B = \emptyset$ and hence
by Lemma \ThisSection.1, $I(G)\cap I(B) \neq\emptyset$.
On the other hand, as $\c{R}$ is
non-trivial, Lemma \ThisSection.4
implies that for each $R\in \c{R}_G$
there is a finite set $F(R)\subset I$ which meets each $I(B)$,
$B\in \c{B}(R)$. Hence,
$$ F \;=\; \cup\{F(R): R\in \c{R}_G\}\cup I(G) $$
is a finite set which satisfies $F\cap I(B)\neq\emptyset$ for
all $B\in \c{B}(R)$ and all $R\in \c{R}$. Recall from
\CPTHM that for the associated {\it finite} product,
$\lambda\Pi_F\mu_i$ contains every open cover of $\Pi_F X_i$.
%If there is a finite subset $\c{R'}\subset \c{R}$ such that
%$\pi_F[\c{B}(\c{R'})]$ contains a finite cover of $\pi_F[X]$,
%then let $\c{R}_X = \c{R'}$;
%otherwise, let $\c{R}_X = \c{R'}$.

Consider the open cover
$$ \c{W}_F \;=\; \pi_F[\c{G}] \wedge \pi_F[\c{B}(\c{R})]$$
of the subproduct $\Pi_F X_i$ of $X$.
(Here
$\c{B}(\c{R})$ denotes the family of all basic sets
contained in some member of $\c{R}$.)
By the result just
mentioned, $\c{W}_F$ belongs to $\lambda(\Pi_F \mu_i)$.
Hence,
there is a Noetherian tree $T_F$
consisting of subsets
of $\Pi_F X_i$ satisfying the uniform
refinement condition with respect to
$\Pi_F\mu_i$ and for which $\End(T_F)\prec \c{W}_F$.
Moreover, we may assume that the elements of $T_F$ are closures of
basic open sets. We extend $T_F$ to a Noetherian tree
$T'_F$ by adding, for each end element $P$, the perverse
product of the trees $T_{\phi_i}(P_i)$ (with respect to
the index set $F$) below $P$, where
$P = \cap\{\pi_i^{-1}[P_i]: i\in F\}$.
(We order the subset $F$ as $\{x_1,\ldots,x_n\}$ and use the
the same perversities as in \CPTHM. The
finite enlargements $F(P)$ defined below will extend
this order. Notice that this limits the level of the
elements from $T_{\phi_i}(P_i)$ to $n$.)
Finally, we
extend $T'_F$ to a tree $T_0$ in $\c{P}(\Pi_I X_i)$ by
crossing each element with $\Pi_{I\setminus F} X_i$.
Then this new tree satisfies the uniform refinement
condition with respect to $\Pi_I\mu_i$.

Let $\End^*T_0$ be the set of all end elements $P$ of $T_0$ for which
there is no $R'\in (\c{R}^{<\omega})^*$, say
$R' = (R_1\cup\cdots\cup R_n)^*$,
with open basic subsets $B\in \c{B}(R)$
such that $B_i\subset R_i$ and $P\subset B_1\cup\cdots\cup B_n$.
If $\End^*T_0\neq\emptyset$, we continue the inductive definition of
$T$. Let $P\in \End^*T_0$. Then $\pi_F[P]$ refines, by the definition
of $T_F$, the cover $\pi_F[\c{G}]$; let $G\in \c{G}$ be such that
$\pi_F[G]\supset \pi_F[P]$.
On the other hand,
$\pi_F[P]$ also refines $\pi_F[\c{W}]$; let $\pi_F[B]\supset\pi_F[P]$,
where $B$ is a basic subset and $B\subset R$ for some $R\in \c{R}$.

Denote by $\c{R}(P)$ the set of
all finite subsets
$\c{R}'\subset\c{R}$, say $\c{R}' = \{R_1,\ldots, R_n\}$, for which
there are $B_i\in \c{B}(\c{R}')$ such that
$$\pi_F[B_1\cup\cdots\cup B_n]\supset \pi_F[P]$$
and additionally $B_i\cap P\neq\emptyset$ for all $i\in \{1,\ldots,n\}$.
Thus, $\c{R}(P)\neq\emptyset$, established by the one-element
subsets.
The set $\c{R}_G$ of all $R\in \c{R}$ which meet $G$ is finite.
Consider an arbitrary element $\c{R}'\in \c{R}(P)$,
$\c{R}' = \{R_1,\ldots,R_n\}$, and a corresponding
set $\{B_1,\ldots, B_n\}$ of basic sets $B_i\in \c{B}(R_i)$
as defined above.
In case
$\c{R}'\not\subset \c{R}_G$,
we have $R_i\cap G=\emptyset$ for some $R_i\in \c{R}'$.
But then
$\pi_F[B_i]\cap \pi_F[G]\neq\emptyset$ and
hence by Lemma \ThisSection.1' we have
$I(B_i)\cap I(G)\not\subset F$. Writing
$\tilde B = B_1\cup\cdots\cup B_n$ and
$I(\tilde B) = I(B_1)\cup\cdots\cup I(B_n)$, we have
$I(\tilde B)\cap I(G)\not\subset F$ for each such
subset $\c{R}'\subset \c{R}_G$ and such
basic sets $B_i\in \c{B}(R_i)$.

On the other hand,
given $\c{R}'\subset \c{R}_G$, Lemma \ThisSection.3'
implies that there is a finite
set $F(\c{R}')\subset I$ such that
$$ (F(\c{R}')\cap I(\tilde B))\cap (I\setminus F)\neq\emptyset,$$
whenever $\tilde B$ and $I(\tilde B)$ are as above.
(Because otherwise we would have $P\subset (\cup \c{R}')^*$,
which is ruled out by the assumption $P\in \End^*(T_0)$.)
Put
$$ F(P) \;=\; \cup\{F(R): \c{R}'\subset \c{R}_G\}\cup I(G)\cup F,$$
and note that $F(P)$ is finite.
Let
$$ \c{W}_P \;=\; \pi_{F(P)}[\c{G}]\wedge \pi_{F(P)}[\c{B}(\c{R})].$$

As before, there is a Noetherian tree $T_P$
consisting of subsets of $P$ satisfying the
uniform refinement condition with respect to
$(\Pi_{F(P)}\mu_i)\upharpoonright P$
and such that $\End(T_P)\prec\c{W}_P$.
For the following and subsequent perverse products,
we order the set $F(P)\supset F$ extending the order
of $F$.
Extend $T_P$ to a Noetherian tree $T'_P$ by adding the
perverse product of the trees $T_{\phi_i}(\pi_i[Q])$
above the elements $Q$, where
$Q\in\End(T_P)$ and $i\in F(P)$.
As above, we enlarge $T'_P$ to a
Noetherian tree in $\c{P}(\Pi_I X_i)$ by
crossing the elements of $T'_P$ with
the product of the $X_i$ for all the remaining factors.
When this is done for each $P\in \End^*(T_0)$,
we obtain a Noetherian tree $T_1$.

In general, if $n > 0$ and $T_n$ has been constructed, then
the inductive step from $T_n$ to $T_{n+1}$ is entirely analogous
to the above construction: We consider the set $\End^*(T_n)$
of all end elements $P$ for which there is no
finite subset $\c{R}'\subset \c{R}$ with
$P\subset (\cup \c{R}')^*$. Given
such an element $Q$, there is a unique $P\in \End^*(T_{n-1})$
below $Q$, and we define $\c{R}(Q)$ as the set of all
finite subsets
$\c{R}'= \{R_1,\ldots, R_n\}\subset \c{R}$ for which there are
basic open sets $B_i\subset R_i$ such that
$\pi_{F(P)}[B_1\cup\cdots\cup B_n]\supset\pi_{F(P)}[Q]$
and $B_i\cap Q\neq\emptyset$ for all $i$.
The construction of
$T'_Q$ then proceeds as above. We set
$$ T_{n+1} \;=\; \cup\{T'_Q: Q\in \End^*(T_n)\}.$$
Finally, we put
$$ T \;=\; \cup\{T_n: n\in \Bbb N\}.$$
We claim that $T$ is Noetherian.
Indeed, suppose to the contrary
that $T$ has an infinite chain $(Q_n: n\in \Bbb N)$. Now each
$T_n$ is Noetherian, so there is in fact an infinite chain
$(P_n: n\in \Bbb N)$ of elements
$P_n\in \End^*(T_n)$ such that
$P_k\supset Q_{n(k)}$. We have $P_{n+1}\subset P_n$
and $F(P_n)\subset F(P_{n+1})$.
Let $J = \cup\{F(P_n): n\in \Bbb N\}$.
The ordering of $T$ along the chain $(P_n)$ is, because of the
perverse orderings of the added trees  $T'_P$,
an ordering for which there is an infinite sequence
$(\pi_i[P_{n(i)}])$ of subsets of $X_i$ for all $i\in J$
such that $\pi_i[P_{(n+1)(i)}]\subset\pi_i[P_{n(i)}]$.
By the construction there are relatively open sets
$U_{n,i} = \phi_i(\pi_i[P_{n(i)}])$ such that
$$ \pi_i[P_{(n+1)(i)}] \subset U_{n,i} \subset \pi_i[P_{n(i)}],$$
where $\phi_i$ denotes the stationary winning strategy
associated with the space $X_i$.

It follows that each
$$ K_i = \bigcap_{n\in \Bbb N} \pi_i[P_{n(i)}]$$ is non-empty
and compact. As $I(P_n)\subset J$ for all $n$, and each $P_n$
is basic, we obtain
$$ \cap_{n\in \Bbb N} P_n \;=\;
    \pi_J[\cap_{n\in \Bbb N} P_n]\times\Pi_{I\setminus J} X_i
       \neq\emptyset$$
and the projection image (denote it by $K$) is compact.

Each $x\in \cap P_n$ has a basic open set $B_x$ such that
$x\in B_x\subset R_x$ for some $R_x\in \c{R}$. There is a
minimal finite set $E$ such that
$$ K = \pi_J[\cap P_n] \subset \cup\{\pi_J[B_x]: x\in E\}.$$
In fact, there is $n$ such that $k\geq n$ implies
$$ \pi_J[P_k]\subset \cup\{\pi_J[B_x]: x\in E\}.$$
(To see this, assume that each $\pi_J[P_k]$ meets the
complement of $\cup\{\pi_J[B_x]: x\in E\}$, say, in a
closed set $A_k$. Then for
each $i\in J$, the sets $\pi_i[A_k]$ form a filter base
finer than $(U_{n,i})$, and thus have cluster points
$x_i\in K_i$. But $K$ is the product of the $K_i$ and
the point $x = (x_i)_{i\in J}$ of $K$ has a neighbourhood
which does not meet any member of the filter base
$(A_k)$, contradicting that $x$ is a cluster point
of the latter.)

Let $\c{R}' = \{R_1,\ldots, R_n\}$,
let $\tilde B = \cup\{B_x: x\in E\}$ and write
$I(\tilde B) = \cup\{I(B_x): x\in E\}$. We have
$\c{R}'\in \c{R}(P_k)$ for all $k\geq n$.

Clearly
$ \pi_{F(P_k)}[P_k] \subset \pi_{F(P_k)}[\tilde B]$
for $k\geq n$.
%Recall that the cover $\c{V}_{P_k}$ of $P_k$
%is defined through
%
%$$ \c{V}_{P_k} \;=\; \pi_{F(P_k)}[\c{G}]\wedge
%      \pi_{F(P_k)}[\c{B}(\c{R})]. $$
%
%The successors of $P_k$ were obtained by considering a
%Noetherian tree $T_{P_k}$ such that
%$\End(T_{P_k})\prec \c{V}_{P_k}$.
%The elements $\pi_{F(P_k)}[R_x]$, $x\in E$,
%form a cover of $\pi_{F(P_k)}[P_k]$, and therefore
%for each end element $P$ of $T_{P_k}$, there is $x\in E$
%such that
%
%$$ \pi_{F(P_k)}[P] \subset \pi_{F(P_k)}[B_x].$$
%
%In particular, there is $x\in E$ with
%
%$$ \pi_{F(P_k)}[P_{k+1}] \subset \pi_{F(P_k)}[B_x].$$
%
On the other hand, $P_{k+1}$ does not refine
$\cup(\c{R}')$ (because it belongs to $\End^*(T_{k+1})$), and we have
$$ (F(\c{R}')\cap I(\tilde B)\setminus F(P_k) \neq\emptyset.$$
But the definition of $F(P_{k+1})$ then implies
$$ (F(P_{k+1})\cap I(\tilde B))\setminus F(P_k)\neq\emptyset.$$
This is valid for all $k\geq n$, and hence
gives the contradiction that
$I(\tilde B)$ is infinite. Thus, $T$ is Noetherian, as claimed.

The end elements $P\in \End(T)$ form a cover of the product
space, and refine the cover $(\c{R}^{<\omega})^*$. The
construction of $T$ implies that the cover
$\End(T)$, and hence $(\c{R}^{<\omega})^*$,
belongs to $\lambda(\Pi_{i\in I}\mu_i)$, and
therefore so does $\c{V}^{<\omega}$. As this is
valid for all normal covers $\c{V}$ of the product space,
we conclude that $\c{V}\in \mu/p\mu$, where
$\mu = \lambda(\Pi_{i\in I}\mu_i)$ and $p\mu$ denotes
the uniformity consisting of all finite covers $\c{U}\in \mu$
(see \referto[Ricea]). But
$\mu/p\mu \subset \lambda\mu = \mu$ and hence
$\c{V}\in \mu$, as desired. $\eop$

\noindent
{\bf Remark \ThisSection.6:} Note that while
the counterpart of \refered[ThNormalCovers] for countably many factors --
\CPTHM -- is stated for regular spaces and pre-uniformities,
\refered[ThNormalCovers] is 'restricted' to uniform spaces. This is needed in the 
last paragraph of the proof, in order to apply \referto[Ricea].
In case we considered directed covers, we could state: Every 
directed normal cover of a product $\Pi X_i$ of regular spaces 
belongs to $\lambda\Pi \mu_i$, where the $\mu_i$ are 
complete pre-uniformities on the spaces $X_i$.

\section{Additional remarks: The locally fine coreflection and `formal topology'.}

One of the motivations for studying `supercompleteness' of
products $\Pi\c{F}X_i$ instead of the mere paracompactness
of the corresponding topological products is the {\it recursive}
or {\it inductive construction} of the refinements of open covers
from the factor covers. This was the motivation behind
the series \referto[Hohtib], \referto[Hohtid], \referto[Hohtie],
\referto[Hohtif], \referto[HohtiPelantb], \referto[HohtiYun]
independently of the questions of `formal spaces'
(\referto[FourmanGrayson], \referto[Sambina]). A similar
motivation can be found in the inductive derivation of
covers for a constructive proof of Tychonoff's Theorem
in the context of frames in \referto[NegriValentini].
Effective presentations of formal spaces have been
studied in \referto[Sigstama].

On the other hand, the `spatiality' (see
\referto[Isbellc] (where it was called `primality'),\referto[Isbelld],
\referto[Plewea]) of the {\it localic products}
$\otimes \Omega(X_i)$ is equivalent to the condition
$$ \lambda(\Pi\c{F}X_i) = \c{F}(\Pi X_i) \leqno(*)$$
whenever the product is paracompact. (Here
$\Omega(X)$ denotes the
topology of $X$ as a locale. The spatiality of the
product $\otimes \Omega(X_i)$ means that it is
isomorphic to the locale $\Omega(\Pi X_i)$ over the
topological product.) Indeed, these three fields have been
brought together in a recent result of the first-named
author \referto[Hohtig] in the following sense: For 
a space $X$, let $\c{O}(X)^*$ denote the (fine) monoid of all 
covers refinable by an open cover. Given a family 
$(X_i)$ of regular spaces, the localic product 
$\otimes\Omega(X_i)$ is spatial if, and only if, 
the equation
$$ \lambda\Pi \c{O}(X_i)^* \; = \; \c{O}(\Pi X_i)^*\leqno(**)$$
holds. This is an analogue of $(*)$, but without any reference
to paracompactness. 

Formal spaces essentially are counterparts of locally fine
pre-uniformities in pre-orders $(P,\leq)$ (partial
order without antisymmetry), defined by `covering relations'
$\mu\subset P\times 2^P$ which satisfy the following axioms:
\item{1)} If $a\in U$, then $\mu(a,U)$.
\item{2)} If $a\leq b$, then $\mu(a,\{b\})$.
\item{3)} If $\mu(a,U)$ and $\mu(a,V)$, then $\mu(a,U\wedge V)$,
          where $U\wedge V$ denotes the set of all $b\in P$
          majorized by both $U$ and $V$.

\item{4)} If $\mu(a,U)$ and for all $u\in U$,
          $\mu(u,V)$, then $\mu(a,V)$.

A set $U\in 2^P$ such that $(a,U)\in \mu$ is considered
a `cover' of the element $a\in P$. By the above axioms,
the covers
of $a$ form a `locally fine' filter for each $a\in P$.
(For a pre-uniformity $\mu$, one takes a filter of
covers of the maximal element (the underlying set)).

Their `effective constructions' are connected -- via the above
theorem -- with the constructions for locally fine coreflections of
products. However, we leave the details of obtaining
effective refinements of open covers of products to future
papers. Let it be mentioned, however, that the definition of
the product of covering relations is reflected in the perverse
product of trees. In the application of perverse products
in Section \refered[FirstMainSection], the immediate successors of a given
element vary with respect to one coordinate only. Furthermore,
the set of perversities is given effectively, and so
is the associated product, relative to the factors.

{\bf Acknowledgement.} The authors are grateful to 
the anonymous referee for valuable remarks, which have 
improved the readability of this paper. 

\bigskip
\SectionBreak
\centerline{\smallcaps References}
\bigskip
{
%\beginref

\ref Alster, K: A class of spaces whose Cartesian product with
     every hereditarily Lindel\"of space is Lindel\"of.-
     Fund.\ Math.\ 114:3, 1981, pp.\ 173--181 [Alstera]

\ref Arhangel'skii, A: On topological spaces complete in the sense
     of \v Cech.- Vestnik Moskov.\ Univ.\ Ser.\ I.\ Mat.\ Mekh.\
     2, 1961, pp.\ 37--40 (Russian) [Arhangelskia]

\ref Arhangel'skii, A: On a class of spaces containing all metric
     and all locally bicompact spaces.- Soviet Math.\ Dokl.\ 4,
     1963, pp.\ 1051 -- 1055 [Arhangelskib]

\ref Dowker, C.\ H., and D.\ Strauss: Sums in the category
     of frames.- Houston J.\ Math. 3, 1977, pp.\ 17--32 [DowkerStrauss]

\ref Fourman, M., and R.\ Grayson: Formal spaces.- The L.\ E.\ J.\
     Brouwer Centenary Symposium, A.\ S.\ Troelstra and D.\ van Dalen
     (eds.), North-Holland, 1982, pp.\ 107--122 [FourmanGrayson]

\ref Friedler, L.\ M., H.\ W.\ Martin and S.\ W.\ Williams:
     Paracompact C-scattered spaces.- Pacific Math.\ Journal 129:2,
     1987, pp.\ 277 -- 296 [Friedlera]

\ref Frol\'\i k, Z: On the topological product of paracompact spaces.-
     Bull.\ Acad.\ Pol.\ Sci.\ Math.\ 8, 1960, pp.\ 747 -- 750
     [Frolika]

\ref Frol\'\i k, Z: Generalizations of $\hbox{\rm G}_{\delta}$-property
     of complete metric spaces.- Czech.\ Math.\ J.\ 10 (85), 1960,
     pp.\ 359--379 [Frolikb]

\ref Ginsburg, S.\ and J.\ R.\ Isbell: Some operators on uniform
     spaces.- Trans.\ Amer.\ Math.\ Soc.\ 93, 1959, pp.\ 145 -- 168
     [Ginsburga]

\ref Hager, A.\ W: Some nearly fine uniform spaces.- Proc.\ London
     Math.\ Soc.\ (3), 28, 1974, pp.\ 517--546 [Hagera]

%
%\ref Heath, R.\ W: A paracompact semi-metric space which is not
%     an $\hbox{\rm M}_3$-space.- Proc.\ Amer.\ Math.\ Soc.\ 17, 1966,
%     pp.\ 868--870.

\ref Hohti, A:On supercomplete uniform spaces.- Proc.\ Amer.\ Math.\ Soc.\ 87:,
     1983, pp.\ 557--560 [Hohtib]

\ref Hohti, A: On Ginsburg-Isbell derivatives and ranks of metric spaces.-
     Pacific J.\ Math.\ 111 (1), 1984, pp.\ 39--48 [Hohtic]

\ref Hohti, A: On supercomplete uniform spaces II.- Czechosl.\ Math.\ J.\ 37,
     1987, pp.\ 376--385 [Hohtid]

\ref Hohti, A: On supercomplete uniform spaces III.- Proc.\ Amer.\ Math. Soc.\
     97:2, pp.\ 339--342 [Hohtie]

\ref Hohti, A: On Relative $\omega$--cardinality and Locally Fine Coreflections
        of Products.- Topology Proceedings 13:1, 1989 [Hohtif]

\ref Hohti, A: Locales, locally fine spaces and formal topology.-
     manuscript, 1998 [Hohtig]

\ref Hohti, A: An Infinitary Version of Sperner's Lemma.- manuscript.
     [Hohtih]

\ref Hohti, A., and Jan Pelant: On complexity of metric
     spaces.- Fund.\ Math.\
     CXXV, 1985, pp.\ 133--142 [HohtiPelanta]

\ref Hohti, A., and Jan Pelant: On supercomplete uniform spaces IV:
     countable products.-  Fund.\ Math. 136:2, 1990, pp.\ 115--120
     [HohtiPelantb]

\ref Hohti, A., and Yun Z: Countable products of \v Cech-scattered supercomplete
     spaces.- Czechosl.\ Math.\ J., to appear [HohtiYun]

\ref Hu\v sek, M., and J.\ Pelant: Extensions and restrictions
     in products of metric spaces.- Topology Appl.\ 25, 1987,
     pp.\ 245 -- 252 [HusekPelant]

\ref Isbell, J: Supercomplete spaces.- Pacific J.\ Math.\ 12, 1962,
     pp.\ 287 -- 290 [Isbella]

\ref Isbell, J: Uniform spaces.- Math.\ Surveys, no.\ 12, Amer.\
     Math.\ Soc., Providence, R.\ I., 1964 [Isbellb]

\ref Isbell, J: Atomless parts of spaces.- Math.\ Scand.\ 31, 1972,
     pp.\ 5--32 [Isbellc]

\ref Isbell, J: Product spaces in locales.- Proc.\ Amer.\ Math.\ Soc.,
     81:1, 1981, pp.\ 116--118 [Isbelld]

\ref Junnila, H., J.\ C.\ Smith, and R.\ Telg\'arsky: Closure-preserving
     covers by small sets.- Topology and Appl.\ 23, 1986, pp.\ 237--262
     [JunnilaSmithTelgarsky]

\ref Kirwan, F: An Introduction to intersection homology theory.-
     Pitman Research Notes in Mathematics Series 187, Longman, 1988
     [Kirwana]

\ref Michael, E: A note on completely metrizable spaces.-
     Proc.\ Amer.\ Math.\ Soc.\ 96:3, 1986, pp.\ 513--522 [Michaelb]

\ref Negri, S., and S.\ Valentini: Tychonoff's theorem in the
     framework of formal topologies.- J.\ Symb.\ Logic, to appear
     in 1997 [NegriValentini]

\ref Pelant, J: Locally fine uniformities and normal covers.-
     Czechosl.\ Math.\ J.\ 37 (112), 1987, pp.\ 181 -- 187
     [Pelanta]

\ref Pelant, J., and M.\ D.\ Rice: On $e$-locally fine spaces.-
     Seminar Uniform Spaces, 1976-77 [PelantRice]

\ref Plewe, T: Localic products of spaces.- Proc.\ London Math.\
     Soc.\ 73:3, 1996, pp.\ 642--678 [Plewea]

\ref Plewe, T:Countable products of absolute $C_\delta$-spaces.-
     Topology Appl.\ 74, 1996, pp.\ 39--50 [Pleweb]

\ref Rice, M.\ D: Metric fine uniform spaces.- J.\ London Math.\
     Soc.\ (2), 11, 1975, pp.\ 53--64 [Riceb]

\ref Rice, M.D: A note on uniform paracompactness.- Proc.\ Amer.\ Math.\
     Soc.\ 62:2, 1977, pp.\ 359--362 [Ricea]

\ref Rudin, M.\ E., and S.\ Watson: Countable products of
     scattered paracompact spaces.- Proc.\ Amer.\ Math.\ Soc.\
     89:3, 1983, pp.\ 551 -- 552 [Rudina]

\ref Sambin, G: Intuitionistic formal spaces and their
     neighbourhood.- Logic Colloquium '88, R.\ Ferro et al
     (eds.), North-Holland, Amsterdam, 1989, pp.\ 261--285 [Sambina]

\ref Sigstam, I: Formal spaces and their effective presentations.-
     Arch.\ Math.\ Logic 34:4, 1995 pp.\ 211 -- 246 [Sigstama]

\ref Stone, A.\ H: Kernel constructions and Borel sets.- Trans.\
     Amer.\ Math.\ Soc.\ 107, 1963, pp.\ 58--70 [Stonea]

\ref Telg\'arsky, R: C-scattered and paracompact spaces.- Fund.\
     Math.\ 73, 1971, pp.\ 59 -- 74 [Telgarskya]

\ref Telg\'arsky, R: Spaces defined by topological games.- Fund.\ Math.\
     LXXXVIII:3, 1975, pp.\ 193--223 [Telgarskyb]

\ref Telg\'arsky, R., and H.\ H.\ Wicke: Complete exhaustive sieves
     and games.- Proc.\ Amer.\ Math.\ Soc.\ 102:1, 1988, pp.\
     737--744 [TelgarskyWicke]

\endref
}

\vskip2cm
\section{The addresses of the authors:}
\bigskip
\hbox{\hfil\vbox{
            \hbox{Aarno Hohti}
            \hbox{University of Helsinki}
            \hbox{Department of Mathematics}
            \hbox{PL 4 (Yliopistonkatu 5)}
            \hbox{00014 Helsingin Yliopisto}
            \hbox{FINLAND}
            \hbox{}}
                   \hskip2cm
                      \vbox{
                         \hbox{Miroslav Hu\v sek}
                         \hbox{Mathematical Institute}
                         \hbox{Charles University}
                         \hbox{Sokolovsk\'a 83}
                         \hbox{186 75 Prague 8}
                         \hbox{Czech Republic}
                         \hbox{}}
                                     \hskip2cm
                                        \vbox{
                                           \hbox{Jan Pelant}
                                           \hbox{Mathematical Institute}
                                           \hbox{The Academy of Sciences of}
                                           \hbox{the Czech Republic}
                                           \hbox{\v Zitn\'a 25}
                                           \hbox{115 67 Prague 1}
                                           \hbox{Czech Republic}}
                                                             \hfil}

\end